\documentclass[fleqn,leqno,keyval,epsfig,11pt,dvips]{amsart}

\usepackage{amssymb}
\usepackage{amsfonts}
\usepackage{amsmath}
\usepackage{amsxtra}
\usepackage{graphicx}

\input epsfx.tex

\newtheorem{theorem}{Theorem}[section]

\newtheorem{proposition}[theorem]{Proposition}
\newtheorem{remark}[theorem]{Remark}
\newtheorem{example}[theorem]{Example}

\title {A New Approach on Constant Angle Surfaces in  $\mathbb{E}^{3}$}
\author[M.I. Munteanu]{Marian Ioan Munteanu}
\author[A.I. Nistor]{Ana Irina Nistor}
\date{}

\address{\scriptsize M. I. Munteanu and A. I. Nistor:
Al.I. Cuza University of Iasi\\
Faculty of Mathematics\\
Bd. Carol I, n. 11\\
700506 - Iasi\\
Romania}

\begin{document}
\maketitle

\begin{abstract}
In this paper we study constant angle surfaces in Euclidean
3--space. Even that the result is a consequence of some classical results
involving the Gauss map (of the surface), we give another approach to
classify all surfaces for which the unit normal makes a
constant angle with a fixed direction. \\[1mm]
{\bf Mathematics Subject Classification (2000):} 53B25\\
{\bf Keywords and Phrases:} constant angle surfaces, Euclidian space.
\end{abstract}

\section{Introduction}

\setcounter{equation}{0}

Recently, constant angle surfaces were studied in product spaces $\mathbb{S}^{2}\times\mathbb{R}$
in \cite{kn:DFVV07} or $\mathbb{H}^{2}\times\mathbb{R}$ in \cite{kn:DM07}, where $\mathbb{S}^{2}$ and $\mathbb{H}^{2}$
represent the unit 2-sphere and the hyperbolic plane, respectively.
The angle was considered between the unit normal of the surface $M$ and the tangent direction to $\mathbb{R}$.
The idea of studying surfaces with different geometric properties in product
spaces was initiated by H. Rosenberg and W. Meeks in \cite{kn:MR04} and \cite{kn:Ros02}, where they have considered the general case
of a surface $\mathbb{M}^{2}$  and they have looked for minimal surfaces properties in the product space
$\mathbb{M}^{2}\times\mathbb{R}$.

In this article we study the problem of constant angle surfaces in Euclidean 3-space. So, we want to find a
classification of all surfaces in Euclidean 3-space for which the unit normal makes a constant angle with a fixed
vector direction being the tangent direction to $\mathbb{R}$.

The applications of constant angle surfaces in the theory of liquid crystals and of layered fluids were
considered by P.Cermelli and A.J. Di Scala in \cite{kn:CS07}, but they used for the study of surfaces another method different
from ours, the Hamilton-Jacobi equation, correlating the surface and the direction field. In \cite{kn:How}, R. Howard explains
how shadow boundaries are formed when the light source is situated at an infinite distance from the surface $M$ using
the geometric model of constant angle surfaces.

\section{Preliminaries}

\setcounter{equation}{0}

Let $\langle\ ,\  \rangle$ be the standard flat metric in $\mathbb{E}^{3}$ and $\widetilde{\nabla}$ its Levi Civita
co\-nnection. We will consider an orientation of $E^{3}$ and denote by $k$ the fixed direction.
Let $M$ be a surface isometrically immersed in $E^{3}$ and denote by $N$ the
unit normal of the surface. Denote by $\theta := \widehat{(N,k)}$, where $\theta \in [0,\pi)$, the
angle function between the unit normal and the fixed direction. A vector is tangent to $M$ if it is
orthogonal to the normal $N$.
\medskip

Recall the Gauss and Weingarten formulas

{\bf (G)}\qquad\qquad $\widetilde\nabla_XY=\nabla_XY+h(X,Y)$

{\bf (W)}\qquad\qquad $\widetilde\nabla_X N=-AX$,

for every $X$ and $Y$ tangent to $M$. Here $\nabla$ is the Levi
Civita connection on $M$, $h$ is a symmetric $(1,2)$-tensor field
taking values in the normal bundle and called the second fundamental
form of $M$ and $A$ is the shape operator. We have
$$
\langle h(X,Y),N\rangle = g(X,AY)
$$
for all $X,Y$ tangent to $M$, where $g$ is the restriction of
the scalar product $\langle\ ,\ \rangle$ to $M$.

\medskip

Decompose $k$ into the tangent and normal part respectively:
\begin{equation}
\label{r1}
\vec{k}=\vec{U}+\cos\theta\vec{N} \ ,\ \rm where\  U\  \rm is\  tangent\  to\  M.
\end{equation}
It follows $\ $
$\|\vec{k}\|^{2}=\|\vec{U}\|^{2}+\cos^{2}\theta\|\vec{N}\|^{2}$ and hence $\|\vec{U}\|=\sin\theta$.

For $\theta \neq 0$, we can define a unit vector field on $M$, namely $e_{1}:= \frac{U}{\|U\|}$.
Let $e_{2}$ be an unitary vector field on $M$ and orthogonal to $e_{1}$. Thus we obtain an orthonormal basis
$\{e_{1}, e_{2}\}$ defined in every point of $M$. From now on we suppose that $\theta$ is constant.

\begin{proposition}
\label{prop:2_1}
In these hypothesis, we have: $[e_{1},e_{2}]\ \|\  e_{2}$.
\end{proposition}
\proof First we calculate $[e_{1},e_{2}]$ and we will notice that it can be written depending only on $e_{2}$.

We will use the following relation:
\begin{equation}
\label{r2}
[e_{1},e_{2}]= \widetilde{\nabla}_{e_{1}}e_{2} - \widetilde{\nabla}_{e_{2}}e_{1}.
\end{equation}
From (\ref{r1}) (and the definition of $e_{1}$) we have $k = \sin\theta e_{1} + \cos\theta N$ and applying
$\widetilde{\nabla}_{e_{2}}$ one gets:
\begin{equation}
\label{r3}
0= \widetilde{\nabla}_{e_{2}}k = \sin \theta \widetilde{\nabla}_{e_{2}}e_{1} + \cos \theta \widetilde{\nabla}_{e_{2}}N.
\end{equation}

Derivating $\langle N , e_{1} \rangle=0$ with respect to ${e_{2}}$ we have the following relation:
\begin{equation}
\label{r4}
\langle \widetilde{\nabla}_{e_{2}}N , e_{1} \rangle + \langle\widetilde{\nabla}_{e_{2}}e_{1} , N \rangle =  0 .
\end{equation}

Weingarten formula yields:
\begin{equation}
\label{r5}
\widetilde{\nabla}_{e_{2}}N = -\rho e_{1} - \lambda e_{2},\rm with \ \rho, \lambda \in C^{\infty}(M).
\end{equation}
From (\ref{r3}) and (\ref{r5}) it follows
\begin{equation}
\label{r6}
\widetilde{\nabla}_{e_{2}}e_{1} = \cot \theta\ (\rho e_{1} + \lambda e_{2}).
\end{equation}

At this point we consider $\theta \neq \frac{\pi}{2}$ (i.e. $\cot\theta \neq 0)$.
The particular case $\theta = \frac{\pi}{2}$ will be treated separately.
\medskip

Combining (\ref{r4}), (\ref{r5}) and (\ref{r6}) we find
$\rho = 0 $ and hence:
\begin{equation}
\label{r7}
\widetilde{\nabla}_{e_{2}}e_{1} = \lambda \cot \theta e_{2} .
\end{equation}

Again, by using the Weingarten formula we have:
\begin{equation}
\label{r8}
\widetilde{\nabla}_{e_{1}}N = -\alpha e_{1}-\gamma e_{2}, \ \rm with \ \alpha, \gamma \in C^{\infty}(M).
\end{equation}
By the same method, applying $\widetilde{\nabla}_{e_{1}}$ to (\ref{r1}) and using (\ref{r8}) we obtain
$$
   \widetilde{\nabla}_{e_{1}}e_{1}=\cot \theta\ (\alpha e_{1} + \gamma e_{2}).
$$
Since $e_{1}$ is unitary it follows that $\alpha$ vanishes.
Moreover, due to the symmetry of the shape operator, i.e.
$\langle Ae_1,e_2\rangle=\langle e_1,Ae_2\rangle$, one immediately gets
that $\gamma$ vanishes too. Hence $Ae_1=0$ and
\begin{equation}
\label{r10}
\widetilde{\nabla}_{e_{1}}e_{1}=0.
\end{equation}
Derivating $\langle e_{1}, e_{2}\rangle = 0$ with respect to $e_1$
and using (\ref{r10}) we get:
\begin{equation}
\label{r11}
\langle \widetilde{\nabla}_{e_{1}}e_{2} , e_{1} \rangle = 0.
\end{equation}
Using the Gauss formula one can write:
$$
   0=\langle Ae_1,e_2\rangle=\langle h(e_1,e_2),N\rangle=
     \langle \widetilde\nabla_{e_1}e_2,N\rangle.
$$
It follows
\begin{equation}
\label{r12}
\widetilde\nabla_{e_1}e_2=0.
\end{equation}
From (\ref{r2}), (\ref{r7}) and (\ref{r12}) we get the following relation for the Lie brackets:
\begin{equation}
\label{r14}
[e_{1} , e_{2}] = -\lambda \cot \theta e_{2} \ ,\ \rm equivalently, \  [e_{1},e_{2}]\ \|\  e_{2}.
\end{equation} \endproof

\rm We conclude this section with the following:

\medskip

\begin{proposition}
\label{prop:2_2}
The Levi Civita connection ${\nabla}$ of $M$ is given by:

\begin{equation}
\label{eq:LC}
\begin{array}{l}
{\nabla}_{e_{1}}e_{1} = 0, \ {\nabla}_{e_{1}}e_{2} = 0,\
{\nabla}_{e_{2}}e_{1} = \lambda \cot \theta\ e_{2},\
{\nabla}_{e_{2}}e_{2} = -\lambda \cot \theta\ e_{1}.
\end{array}
\end{equation}
\end{proposition}
\proof \rm The expressions can be obtained by straightforward computations.
See also \cite{kn:DFVV07} and \cite{kn:DM07}.
\endproof

\section{The characterization of constant angle surfaces}

\setcounter{equation}{0}

Due to Proposition \ref{prop:2_1} one can choose now a local coordinate system in each point of the surface
$M$, namely a parametrization:
$$r = r(u,v)= (x(u,v), y(u,v), z(u,v))$$
 such that the tangent vectors are:
$r_{u} = e_{1}$ and $r_{v} \| e_{2}$. Let $r_{v}:= \beta (u,v)e_{2}$, where $\beta$ is a smooth function on $M$.
Hence, the metric on $M$ can be written as:
\begin{equation}
\label{r16}
g = du^{2} + \beta^{2}(u,v)dv^{2} .
\end{equation}

\begin{remark}
\label{prop:3_1}
 The coefficients of the first fundamental form are: $E = 1$, $F = 0$,
$G = \beta^{2}(u,v)$.
\end{remark}

\rm From Proposition \ref{prop:2_2} one can write now the Levi Civita connection of $M$ in terms of the coordinates $u$ and
$v$. It follows that the parametrization $r$ satisfies the following PDE's:
\begin{equation}
\label{r17}
r_{uu} = 0
\end{equation}
\begin{equation}
\label{r18}
r_{uv} = \frac{\beta_{u}}{\beta}\ r_{v}
\end{equation}
where $\beta$ satisfies the following PDE :
\begin{equation}
\label{r19}
  \beta_{u} - \beta \lambda \cot \theta = 0
\end{equation}
and finally,
\begin{equation}
\label{r20}
r_{vv} = \frac{\beta_{v}}{\beta}\ r_{v} + \beta^{2}\lambda \cot \theta r_{u} - \beta^{2} \lambda N.
\end{equation}



Using the Schwartz identity: $\widetilde{\nabla}_{\partial_{u}}\widetilde{\nabla}_{\partial_{v}}N =
\widetilde{\nabla}_{\partial_{v}}\widetilde{\nabla}_{\partial_{u}}N $ and the ex\-pre\-ssions of the partial
derivatives of the unit normal of the surface $M$: $N_{u}=0$ and $N_{v}=-\lambda r_{v}$, we have that $\lambda$
satisfies the following PDE:
\begin{equation}
\label{r21}
\lambda_{u}+\lambda^{2}\cot \theta = 0 .
\end{equation}

Now we have to find the functions $\lambda$ and $\beta$ in order to write the parametrization $r$ of the surface
$M$.

\begin{remark}
\label{rem:3_1}
Since $N_{u}=0$ it follows that the coefficients of the se\-cond fundamental form: $e = f = 0$. This fact
implies that the Gaussian curvature  of $M$ vanishes, $K = 0$. So, the surface $M$ is locally flat.
\end{remark}
\begin{remark}
In terms of the Gauss map of the surface we can say that it makes a constant angle with a fixed direction,
which is equivalent to the fact that the Gauss map lies on a circle in the sphere $S^2$. Since it has no interior
points in $S^2$, it follows that the Gaussian curvature of the surfaces vanishes identically.
\end{remark}
\begin{proposition}
\label{prop:3_2}
The functions $\lambda$ and $\beta$ are given by the following expressions:
\begin{equation}
\label{r24}
\lambda(u,v) = \frac{\tan \theta}{u + \alpha(v)}
\end{equation}
\begin{equation}
\label{r25}
\beta(u,v) = \varphi(v)(u+\alpha(v)),
\end{equation}
where $\alpha$ and $\varphi$ are smooth functions on $M$
or,
\begin{equation}
\label{r22}
\lambda(u,v) = 0
\end{equation}
\begin{equation}
\label{r23}
\beta(u,v)=\beta(v).
\end{equation}
\end{proposition}
\proof
First we solve (\ref{r21}) and we find $\lambda$ and then we substitute it
in (\ref{r19}) obtaining $\beta$.
\endproof
\begin{theorem}
\label{th:3_1}
{\rm (of characterization)}
\it A surface $M$ in $\mathbb{E}^{3}$ is a constant angle surface if and only if it is locally isometric
to one of the following surfaces:\\[2mm]
$(i)$ either a surface given by
\begin{equation}
\label{r26}
r:M \rightarrow \mathbb{E}^{3},\ (u,v)\mapsto (u\cos\theta(\cos v,\sin v)+\gamma(v), u\sin\theta)
\end{equation}
with
\begin{equation}
\label{eq:gamma_case1}
\gamma(v)=\cos\theta\ \big(-\int\limits_0^v \alpha(\tau)\sin\tau d\tau,\int\limits_0^v\alpha(\tau)\cos\tau d\tau\big)
\end{equation}
for $\alpha$ a smooth function on an interval $I$,\\[1mm]
$(ii)$ or an open part of the plane $x\sin\theta-z\cos\theta=0$,\\[1mm]
$(iii)$ or an open part of the cylinder $\gamma\times\mathbb{R}$, where $\gamma$ is a smooth curve in $\mathbb{R}^2$.\\[1mm]
Here $\theta$ is a real constant.
\end{theorem}
\proof First we prove that all these surfaces define indeed a constant angle surface in
$\mathbb{E}^{3}$. Item $(ii)$ is obvious and item $(iii)$ corresponds to $\theta=\frac\pi2$.
For item $(i)$ we have the tangent vectors
$$
\begin{array}{l}
r_u=\big(\cos\theta\cos v,\cos\theta\sin v,\sin\theta\big)\\[1mm]
r_v=\big((u+\alpha(v))\cos\theta(-\sin v,\cos v),0\big).
\end{array}
$$
Thus, the unit normal is $N=\big(-\sin\theta(\cos v,\sin v),\cos\theta\big)$ and hence,
the angle between $N$ and and the fixed direction $k$ is the constant $\theta$.

\medskip

Conversely, we have to prove that a constant angle surface in $\mathbb{E}^{3}$ is
as in the statement of the theorem.
Since $e_1=r_u$, from (\ref{r1}) we get
$$
k=\sin\theta\ r_{u} + \cos\theta\ N.
$$
Using Remark \ref{prop:3_1} and from the previous relation it easily follows that
 $\langle r_{u},k\rangle=\sin\theta$ and $\langle r_v,k\rangle=0$. Hence
the third component of $r(u,v)$ is $z(u,v)=u\sin\theta$.
\medskip

At this point, the parametrization of $M$ becomes:
\begin{equation}
\label{r29}
r(u,v)=(h(u,v),u\sin\theta)
\end{equation}
where $h(u,v)=(x(u,v),y(u,v)) \in \mathbb{R}^{\rm 2}$.\\

We analyze the two cases for $\lambda$ and $\beta$ furnished by the Proposition \ref{prop:3_2}.\\

{\bf CASE I.}

Since $r_{uu}=0$ we have $h_{uu}=0$. On the other hand, $e_1=r_u=(h_u,\sin\theta)$ is a unit vector,
which means that $|h_u|=\cos\theta$. Hence $h_u=\cos\theta f(v)$, where $f(v)\in\mathbb{R}^2$ and
$|f(v)|=1$ for any $v$, i.e. $f$ is a parametrization of the circle $S^1$. By integration we obtain
$$
    h(u,v)=u\cos\theta f(v)+\gamma(v)
$$
where $\gamma$ is a smooth curve in $\mathbb{R}^2$.\\[1mm]
It follows that $r_v=(u\cos\theta f'(v)+\gamma'(v),0)$. Since $r_{uv}=\frac{\beta_u}\beta\ r_v$
we get $\gamma'(v)=\cos\theta \alpha(v) f'(v)$.\\[1mm]
Without loss of the generality we can suppose that $f$ is the natural
parametrization for $S^1$, i.e. $f(v)=(\cos v,\sin v)$ (this corresponds to a change of the parameter $v$).\\[2mm]
One obtains the parametrization for $M$
$$
r(u,v)=\big(u\cos\theta(\cos v,\sin v)+\gamma(v),u\sin\theta\big)
$$
where $\gamma$ is given by (\ref{eq:gamma_case1}).

\medskip

{\bf Case II.}

Due to $r_{uu}=0$ and $r_{uv}=0$ it follows that $h_{uu}=0$ and $h_{uv}=0$, which imply that $h_u$ is a constant vector in
$\mathbb{R}^2$ of length $\cos \theta$, i.e. $h_u=\cos\theta(\cos\mu,\sin\mu)$, $\mu\in\mathbb{R}$. Hence
$$
   h(u,v)=u\cos\theta(\cos\mu,\sin\mu)+\gamma(v)
$$
where $\gamma$ is a smooth curve in $\mathbb{R}^2$.

Recall that $r_u$ and $r_v$ are orthogonal. Consequently,
$$\gamma(v)=\alpha(v)(-\sin\mu,\cos\mu), \ \alpha\in C^\infty(I).\qquad\qquad (*)$$
The parametrization of $M$ can be written as
$$
r(u,v)=\big( u\cos\theta(\cos\mu,\sin\mu)+\gamma(v),u\sin\theta\big)
$$
with $\gamma$ given from $(*)$.\\[2mm]
A rotation of angle $\mu$ in the plane $(x,y)$ yields the following parametrization for $M$
\begin{equation}
\label{eq:param_case2}
r(u,v)=(u\cos\theta, \alpha(v),u\sin\theta)
\end{equation}
which parameterizes the plane $x\sin\theta-z\cos\theta=0$.\\

{\bf Particular cases for the constant angle $\theta$:}\\

$\bullet$ $\theta = 0$ : the normal $N$ coincides with the direction $k$.
Since $r_u$ and $r_v$ are tangent to $M$ it follows $\langle r_u,k\rangle=0$
and $\langle r_v,k\rangle=0$ and thus $\langle r,k\rangle={\rm constant}$.
This is the equation of a plane parallel to $(x,y)-$plane.
It can be parameterized as $r(u,v)=(u,v,0)$.\\

$\bullet$ $\theta = \frac{\pi}{2}$ : $k$ is tangent to the surface.
In this case $M$ is the product of a curve in $\mathbb{R}^{2}$ and $\mathbb{R}$ (cylindrical surface),
which can be parameterized as in (\ref{r26}) by:
$r(u,v)= (\gamma(v), u)$, where $\gamma(v) \in \mathbb{R}^{2}$.

Now the  theorem is completely proved.
\endproof

We give some examples of constant angle surfaces, parameterized by
(\ref{r26}) for different functions $\alpha$ in (\ref{eq:gamma_case1}).
All pictures are realized by using Matlab.

\begin{example}\rm
In all the following four examples we consider $\theta=\frac\pi4$.\\[2mm]
1. $\alpha(v)=1$:

\quad $r(u,v)=\frac1{\sqrt{2}}\ \big((1+u)\cos v-1, (1+u)\sin v, u\big)$\\[1mm]

2. $\alpha(v)=v$:

\quad $r(u,v)=\frac1{\sqrt{2}}\ \big((u+v)\cos v-\sin v,(u+v)\sin v+\cos v-1, u\big)$


\begin{figure}[htb]
\begin{center}
\epsfxsize=60mm
\centerline{\leavevmode
\epsffile{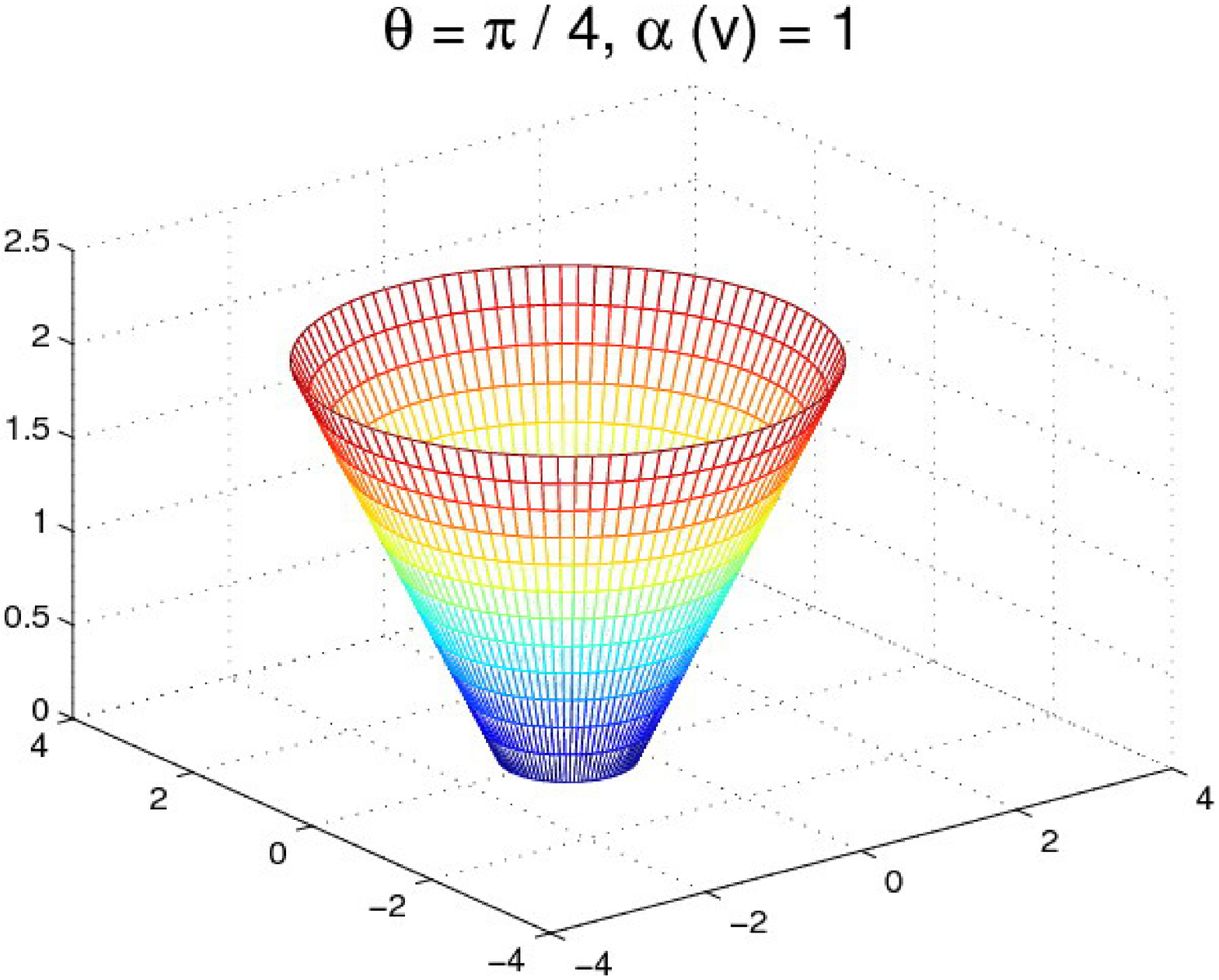}
\ \
\epsfxsize=60mm
\epsffile{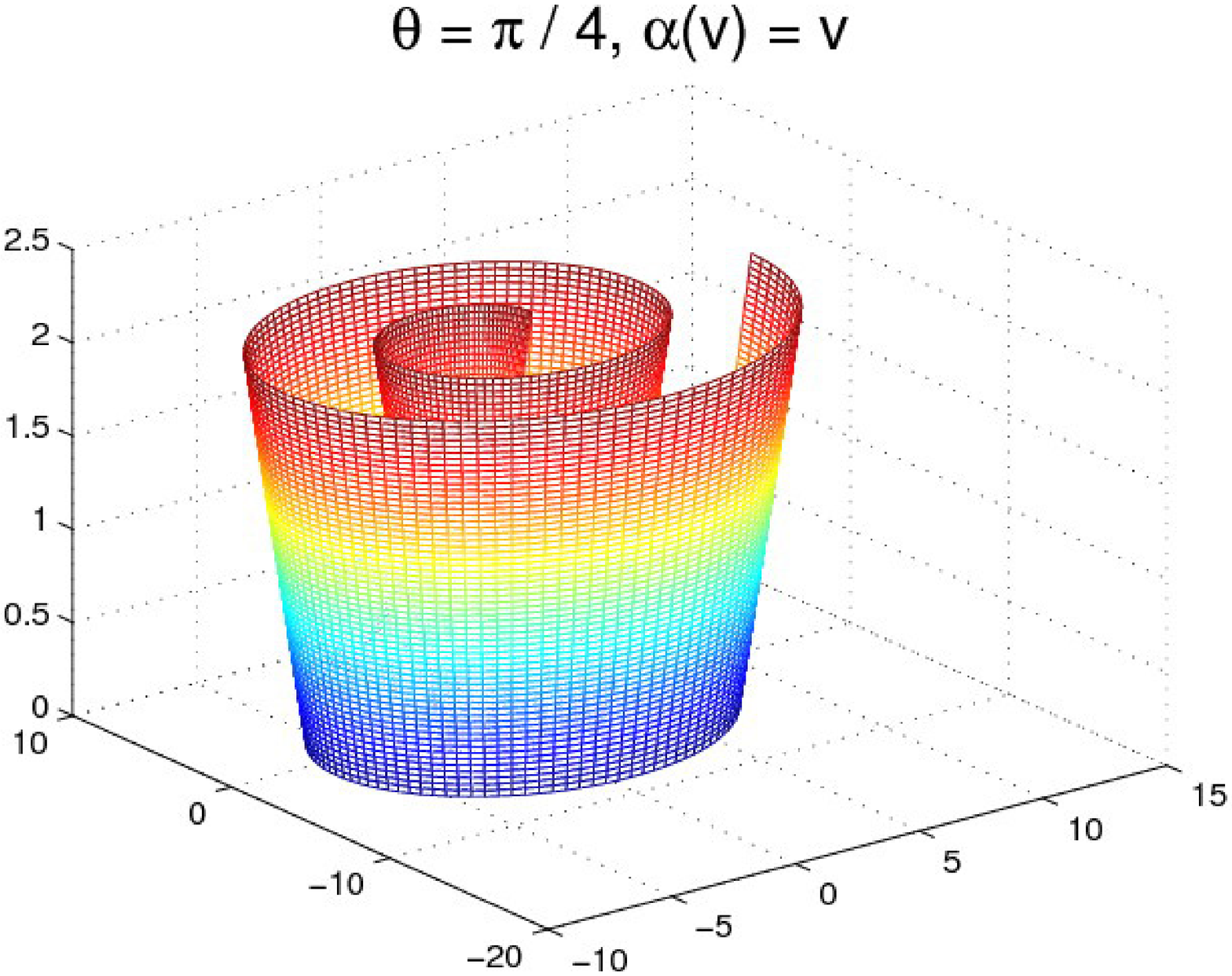}
}
\end{center}
\caption{}
\end{figure}

3. $\alpha(v)=\cos v$:

\quad $r(u,v)=\frac1{\sqrt{2}}\ \big(u\cos v-\frac{\sin^2v}2,u\sin v+\frac{v+\sin v\cos v}2, u\big)$\\[1mm]

4. $\alpha(v)=2\sin v$:

\quad $r(u,v)=\frac1{\sqrt{2}}\ \big(u\cos v-v+\cos v\sin v,u\sin v+\sin^2 v, u\big)$\\[2mm]

\begin{figure}[htb]
\begin{center}
\epsfxsize=60mm
\centerline{\leavevmode
\epsffile{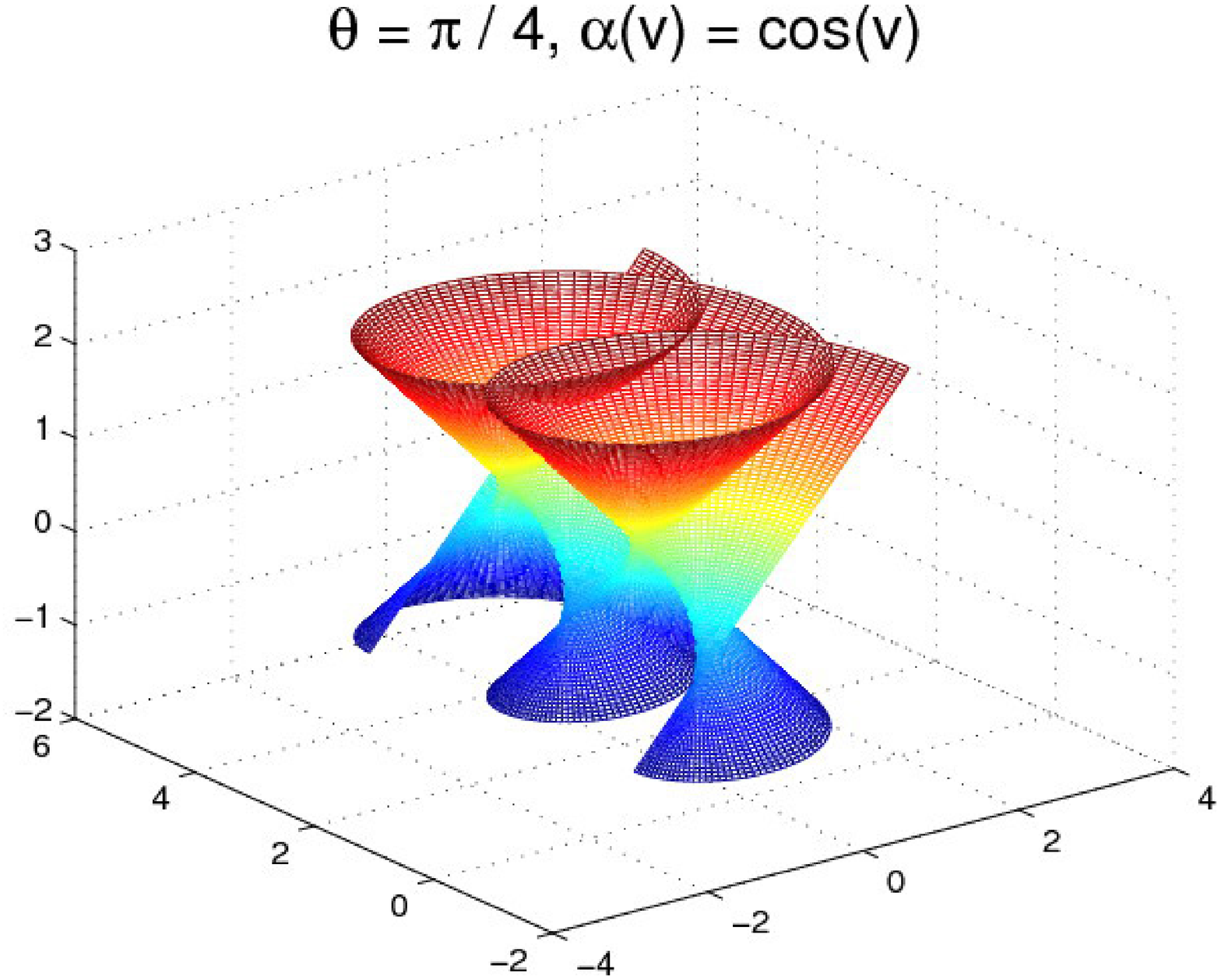}
\qquad
\epsfxsize=60mm
\epsffile{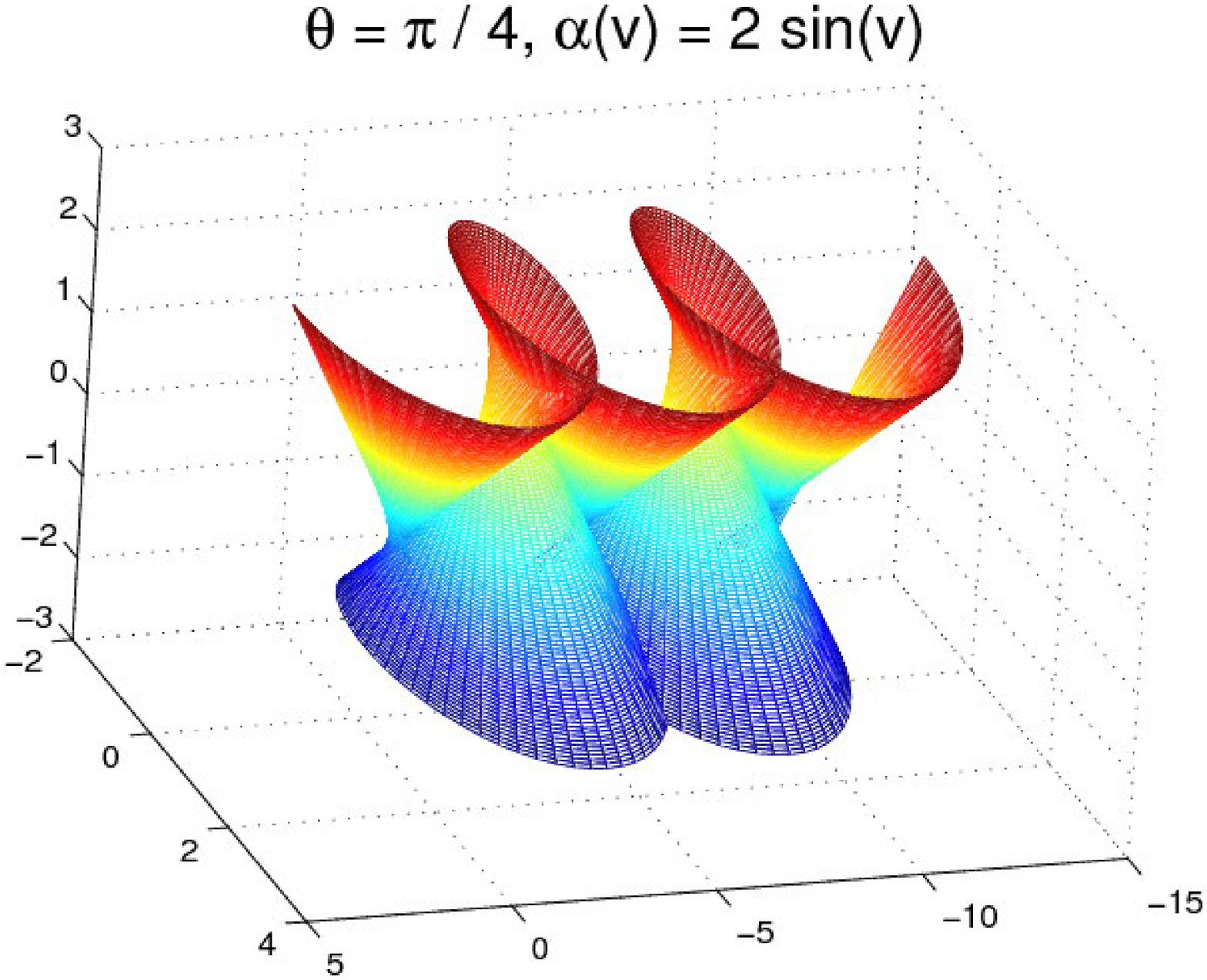}
}
\end{center}
\caption{}
\end{figure}
\end{example}

\medskip

We give now the following result:

\begin{proposition}
\label{prop:3_3}\

{\bf 1.} The only minimal constant angle surfaces in Euclidean $3$-space are the planes which make the angle
$\theta$ with the  fixed direction $k$.\\[2mm]
{\bf 2.} The constant angle surfaces in Euclidean $3$-space with non-zero constant mean curvature are the cylindrical
surfaces.
\end{proposition}
\proof Recall the formula $H=\frac{1}{2}\ \frac{eG-2fF+gE}{EG-F^{2}}$.
Using Remark \ref{prop:3_1} and Remark \ref{rem:3_1} we get that
\begin{equation}
\label{r39}
H=\frac{1}{2}\frac{g}{\beta^{2}(u,v)}\ .
\end{equation}
Looking for all minimal surfaces (i.e. $H=0$) we should have $g=0$. Now we refind here the case $\lambda = 0$
corresponding to the planes which make the angle $\theta$ with the  fixed direction $k$.

For the second statement ($M$ is CMC), (\ref{r39}) implies that $\lambda={\rm constant}$. But $\lambda$ satisfies (\ref{r21})
so we must have $\theta = \frac{\pi}{2}$. In this particular case we found the cylindrical surfaces $\gamma\times\mathbb{R}$,
$\gamma$ smooth curve in $\mathbb{R}^2$.
\endproof

\section{Conclusions}

We can compare now all three results obtained for different ambient spaces, namely for
$\mathbb{S}^{2}\times\mathbb{R}$, $\mathbb{H}^{2}\times\mathbb{R}$ and $\mathbb{E}^{3}$, respectively.
Thus we have: {\em $M$ is a constant angle surface if and only if it is given by an immersion $r$
of the following form}

1. $r:M \rightarrow \mathbb{S}^2 \times \mathbb{R}$,

\quad $(u,v) \mapsto (\cos(u \cos\theta)f(v)+\sin(u \cos\theta)f(v)\times f'(v),\ u \sin\theta)$

\quad where $f:I\rightarrow\mathbb{S}^{2}$ is an unit speed curve in $\mathbb{S}^{2}$ - the unit 2-sphere

\quad and $"\times"$ is the vector cross product in $\mathbb{R}^3$;

2. $r : M \rightarrow \mathcal{H} \times \mathbb{R}$,

\quad $(u,v) \mapsto (\cosh(u \cos\theta)f(v)+\sinh(u \cos\theta)f(v)\boxtimes f'(v),\ u \sin\theta)$

\quad where $f:I \rightarrow\mathcal{H}$ is an unit speed curve on $\mathcal{H}$ - the hyperboloid model of

\quad $\mathbb{H}^{2}$
and $"\boxtimes"$ is the Lorentzian cross product in $\mathbb{R}^3_1$ - Lorentzian 3-space;

3. $r : M \rightarrow \mathbb{E}^{3}$,\
$(u,v) \mapsto (u \cos\theta f(v)+\gamma(v),\ u\sin\theta)$

\quad
where $f:I \rightarrow\mathbb{R}^{2}$ is a parametrization of the unit circle $\mathbb{S}^1$, or $f$ is
a unit

\quad constant vector and $\gamma'(v)\perp f(v)$.

\begin{remark}\rm
\label{rem:4_1}
The third component (along $\mathbb{R}$) in all of these cases is the same: $ z(u,v)= u\sin\theta$.
\end{remark}
\begin{remark}\rm
\label{rem:4_2}
In $\mathbb{S}^2 \times \mathbb{R}$ the surface $M$ has the constant Gaussian curvature
$K=\cos^{2}\theta > 0$, in $\mathcal{H} \times \mathbb{R}$ one gets $K=-\cos^{2}\theta < 0$ while in $\mathbb{E}^{3}$
it vanishes ($K=0$).
\end{remark}

\section{Appendix}

{\bf Applications to the theory of liquid crystals.}
In terms of differential geometry, we studied the constant angle surfaces in $E^{3}$ whose unit
normal forms a constant angle with an assigned direction field. From the point of view of physics, this geometric
condition is equivalent to an Hamilton-Jacobi equation correlating the surface and the direction field.

In the physics of
interfaces in liquid crystals and of layered fluids, these surfaces are studied when the direction field, in our
case
$k$, is singular along a line or a point. We can see in \cite{kn:CS07} how constant angle surfaces may be used to describe
interfaces
occurring in special equilibrium configurations of {\em nematic}
and {\em smectic C liquid crystals},
and to determine the shape of disclination cores in nematics. The last aspect, applications of constant angle surfaces
in nematics was developed by E.G. Virga in \cite{kn:Vir89}, and more recently, for example,
by P. Prinsen and P. van der Schoot in \cite{kn:PS03}, \cite{kn:PS04}, \cite{kn:PS04'}.

\medskip

{\bf Acknowledgement.}
The first author was supported by grant CEEX-ET n. 69/5871/2006-2008, ANCS, Romania.

\end{document}